\documentclass{amsart}

\usepackage{amsmath,amssymb,latexsym, amscd}
\usepackage{exscale, cite, epsfig, graphics, hyperref}

\renewcommand{\geq}{\geqslant}
\renewcommand{\leq}{\leqslant}

\newtheorem{theorem}{Theorem}[section]

\newtheorem{example}[theorem]{Example}
\newtheorem*{main-theorem}{Main Theorem}
\newtheorem*{remark*}{Remark}
\numberwithin{equation}{section}

\newcommand{\be}{\begin{equation}}
\newcommand{\ee}{\end{equation}}
\newcommand{\bse}{\begin{subequations}}
\newcommand{\ese}{\end{subequations}}

\usepackage{hyperref,xcolor}

\title[Pressure transfer functions]
{Pressure transfer functions for\\interfacial fluids problems}

\author[Chen]{Robin~Ming~Chen}
\address{Department of Mathematics, University of Pittsburgh, Pittsburgh, PA 15260 USA} 
\email{mingchen@pitt.edu}  

\author[Hur]{Vera~Mikyoung~Hur}
\address{Department of Mathematics, University of Illinois at Urbana-Champaign, Urbana, IL 61801 USA}
\email{verahur@math.uiuc.edu}

\author[Walsh]{Samuel~Walsh}
\address{Department of Mathematics, University of Missouri, Columbia, MO 65211 USA} 
\email{walshsa@missouri.edu}  

\keywords{Water waves; vorticity; stratification; pressure; transfer function}
\subjclass[2010]{35R35,35C07,76B15,76B47,76B55,76B70}

\begin{document}

\maketitle

\begin{abstract}
We make a consistent derivation, from the governing equations, 
of the pressure transfer function in the small-amplitude Stokes wave regime 
and the hydrostatic approximation in the small-amplitude solitary water wave regime,
in the presence of a background shear flow. 
The results agree with the well-known formulae in the zero vorticity case,
but they incorporate the effects of vorticity through solutions to the Rayleigh equation.
We extend the results to permit continuous density stratification
and to internal waves between two constant-density fluids. Several examples are discussed. 
\end{abstract}

\section{Introduction}\label{sec:intro}
We study the reconstruction from pressure problem for traveling water waves with vorticity or density stratification. 
Consider a body of water whose motion is governed by the steady Euler equations for an incompressible fluid in two dimensions (see Section~\ref{sec:vorticity}).  Assume that the fluid region is bounded above by a free surface  
and below by a fixed horizontal bed.
We then ask:  Can one determine the surface profile and the velocity field inside the fluid knowing only the trace of its pressure at the bed and the background current?

The question is intriguing from a mathematical perspective, but it originates from practical applications in ocean engineering.  The vastness of the ocean seriously limits the types and amounts of data that can realistically be gathered.  Tracking surface waves thus requires a reconstruction method that uses minimal information and, in particular, one that can be employed in rough seas or over large areas.  An attractive possibility is to use pressure transducers, which sample the pressure near the sea bed at key positions in the ocean.  Indeed, they have emerged as one of the primary tools for tsunami monitoring.  

Remarkably, reconstruction from pressure is mathematically possible, at least in some physical regimes.  For instance, it was proved in \cite{henry2013pressure} that 
once the wave speed and the vorticity-stream function relation are specified, 
the trace of the pressure at the bed uniquely determines the surface wave. 
In \cite{ChenWalsh2015pressure}, two of the  authors extended the result 
to permit stratified fluids and relaxed the regularity assumptions.  
The idea behind these works is that, from Bernoulli's law restricted to the bed, 
the trace of the pressure furnishes Cauchy data for an elliptic PDE 
whose solution describes the motion in the fluid region. 
If a solution is known to exist --- as we must assume in order to begin a reconstruction procedure --- 
then we conclude that there is at most one surface wave for the prescribed pressure data.  
We call the pressure-to-surface mapping the (nonlinear) \emph{pressure transfer function}.  
In the internal wave regime, likewise, the interfaces between adjacent layers are 
each uniquely determined by the trace of the pressure at the bed. 
We also refer to this pressure-to-interface(s) mapping as the pressure transfer function. 

When the flow in the bulk is irrotational, that is, the velocity field is curl free,
the elliptic PDE reduces to the Laplace equation in a strip, which can then be solved explicitly. 
This was observed in \cite{Constantin2012}, for instance, where an exact relation was derived 
between the trace of the pressure at the bed and the surface of a solitary water wave. 
The result was extended in \cite{clamond2013recovery,clamond2013new} to periodic waves. 
Based on an alternative formulation of the problem, 
another exact relation was derived in \cite{OVDH2012} for solitary and periodic water waves.
Two of the authors in \cite{ChenWalsh2015pressure} obtained similar results for layer-wise irrotational flows.

The main appeal of these works is that 
they are concerned with exact solutions of the problem without asymptotic assumptions.  
In particular, the formulae apply to large-ampltiude waves.   Though they are typically implicit, there are straightforward numerical implementations, 
and the results to varying degrees match laboratory experiments; see, e.g., \cite{deconinck2012relating}.  

However, there are serious barriers to \emph{proving} that the numerical schemes are convergent.  
While analytically the pressure at the bed uniquely determines the surface wave, 
the correspondence is not a continuous mapping between conventional function spaces. 
One can think of this as a consequence of the ill-posedness of the Cauchy problem for an elliptic PDE.
Indeed, the exact formula in \cite{Constantin2012}, for instance, is only valid in the strong sense
when the trace of the pressure at the bed has Fourier coefficients that decay exponentially fast. 
This is true for an actual solution in the zero vorticity case (cf., e.g., \cite{PlotnikovToland2002}), 
but the set of functions enjoying this property is not open in any natural topology.  
Therefore, when one introduces noise or makes approximations necessary to do computations, 
there is no mathematical justification for convergence. 
Moreover obtaining a posteriori bounds, or otherwise quantifying the error, appears to be highly intractable.

Another issue is that, at present, there are no exact expressions 
relating the trace of the pressure at the bed and the surface wave in rotational flows and stratified fluids
even though we know that pressure transfer functions exist in these regimes.  

To summarize, the exact reconstruction methods work well analytically, 
and they appear to predict reasonably well in practice, 
but it is very difficult to prove that the numerical schemes are convergent. 
Moreover they do not apply as widely as one might hope.  
With this in mind, the present objective is to derive linear transfer functions 
for various classes of surface and interfacial waves, which currently do not have nonlinear analogues. 
To the best of the authors' knowledge, these are new. 
We exploit the rigorous existence theory, rather than relying on a formal approximation of the problem 
(see, e.g.,  \cite{ES2008}). 
We thereby sacrifice the ability to handle large-ampltiude waves,
but in return gain precise bounds on the error committed by the linear approximation.  
Moreover one may generalize the method to gain improved accuracy 
by using higher order perturbations of the governing equations (see, e.g., \cite{HL2015pressure}) 
or higher order expansions of the solution, which may be computed from the existence theory 
(see, e.g., \cite{constantin2015approximations}).  

Linear transfer functions enjoy a long history in the literature; 
see, e.g., \cite{bishop1987measuring, ES2008, kinsman1965wind, kuo1994transfer, schiereck2003introduction}.  
Our contribution here is to allow for general distribution of vorticity
coming, for instance, from a background current or density stratification.  
This substantially complicates the reconstruction from pressure problem.  
Indeed, the pressure distribution in a water wave with vorticity may differ 
in interesting and counterintuitive ways from the zero vorticity case.  
For instance, it was proved in \cite{constantin2010pressure} that 
the pressure of a Stokes wave in the irrotational flow is maximized at the bed at the point directly below the crest 
and decreases monotonically as one moves upward, or to the right or left.  
However, numerical experiments in \cite{da1988steep,ko2008large,ko2008vorticity,vasan2014pressure} suggest that this is not true for some vorticity distributions.  Our results show that linear transfer functions incorporate the effects of vorticity in a subtle way, involving solutions of the Rayleigh equation. 

In Section~\ref{sec:vorticity} we introduce the governing equations 
for traveling water waves in two dimensions with vorticity, 
and we derive the pressure transfer functions for small-amplitude Stokes waves 
and the hydrostatic approximation for small-amplitude solitary water waves. 
In Section~\ref{sec:continuous} and Section~\ref{sec:internal} we extend the results 
to surface waves in continuously stratified fluids and to internal waves in layerwise-constant density fluids. 
Throughout we make an effort to present concrete examples.  

\section{Surface water waves with vorticity}\label{sec:vorticity}
Consider an incompressible inviscid fluid in two dimensions lying below a body of air and subjected to the force of gravity.  Our interest lies in determining the wave motion at the interface between the air and fluid. 
Adopting Cartesian coordinates $(x,y)$, we take the $x$-axis pointing in the direction of wave propagation and the $y$-axis vertically upward. The fluid at time $t$ occupies the connected region $\Omega(t)\subset\mathbb{R}^2$, 
which is bounded above by a {\em free} surface $S(t)$ and below by a fixed horizontal bottom $\{y=0\}$. 
We assume that the free surface has graph geometry, that is, 
\[
S(t)=\{(x,y)\in\mathbb{R}^2:y =h(x;t)\}\qquad\textrm{for all }t > 0,
\]
for some single-valued, non-negative, and smooth surface profile $h$. 

Let $(u,v) = (u(x,y;t),v(x,y;t))$ be the velocity field and $P = P(x,y;t)$ be the pressure. They satisfy the Euler equations for an incompressible fluid:
\begin{equation}\label{E:euler1} \left\{
\begin{split}
&u_t+uu_x+vu_y=-\frac{1}{\rho}P_x,\\
&v_t+uv_x+vv_y=-\frac{1}{\rho}P_y-g\qquad\text{in }\Omega(t),\\
&u_x+v_y=0,
\end{split} \right.\end{equation}
where $g$ is the constant acceleration of gravity. Although an incompressible fluid such as water may have variable density, 
in many applications it is sufficient to assume that the density is constant.  In this section, we shall take $\rho = 1$. The case of stratified fluids will be considered in Section~\ref{sec:continuous} and Section~\ref{sec:internal}.

The kinematic boundary conditions
\begin{equation}\label{E:top-bdry1}
v=h_t+uh_x\qquad\text{on }S(t)
\end{equation}
and 
\begin{equation}\label{E:b-bdry1}
v=0\qquad\text{on }\{y=0\}
\end{equation}
ensure that fluid particles on the boundaries are confined to them at all times. The dynamic boundary condition
\begin{equation}\label{E:top-bdry2}
P=P_{\textrm{atm}}\qquad\text{on }S(t)
\end{equation}
states that the pressure is continuous across the surface, and hence it agrees with the constant atmospheric pressure there. This neglects the effects of surface tension and the dynamics in the air.

\subsection{Introducing a shear flow and scaling the variables}\label{ss:scaling}

The local spin of a fluid element is measured by the vorticity of the flow, which in two dimensions is defined by
\begin{equation}\label{def:vorticity}
\omega=v_x-u_y.
\end{equation}
Since $\Omega(t)$ is simply connected, we may write that
\begin{equation}\label{def:psi}
\psi_x=-v\quad\textrm{and}\quad\psi_y=u,
\end{equation}
for a stream function $\psi$. The fluid is said to be irrotational provided that $\omega \equiv 0$.  Because a flow that is initially irrotational remains so for all time (see \cite{Johnson1997}, for example), it is a common practice to suppose that $\omega$ vanishes identically.  
This choice is more attractive in view of the fact that irrotational flows are much easier to handle --- 
both analytically and numerically --- than rotational flows.  
But the effects of vorticity are physically significant in many circumstances. Vorticity arises, for instance, from density stratification, the action of wind blowing over the water, an incoming sheared current, or the effects of seabed topography. 
The study of rotational flows necessitates a careful consideration of 
the dynamics in the bulk of the fluid and its intimate coupling to the motion of the boundary. 

Note that \eqref{E:euler1}--\eqref{E:top-bdry2} admits a solution of the form
\begin{equation}\label{def:shear}
h\equiv h_0,\quad(u,v)=(U(y),0)\quad\textrm{and}\quad P =P_{\textrm{atm}}+g(h_0-y)
\end{equation}
for arbitrary $h_0>0$ and $U$ in the $C^1\in([0,h_0])$.  Physically, this represents a shear flow, for which 
the velocity and the free surface are horizontal and the pressure is hydrostatic. The vorticity generated by \eqref{def:shear} is $ -U'(y)$.
Here we are interested in waves propagating in the $x$-direction over a prescribed shear flow of the form. 
In the remainder of the section, therefore, $h_0$ and $U$ are held fixed.


In order to systematically characterize various types of approximations, we work with two dimensionless parameters: 
\begin{equation}\label{def:parameters}
\delta=\textrm{the long wavelength parameter}\quad\textrm{and}\quad
\epsilon=\textrm{the amplitude parameter},
\end{equation}
and define a set of rescaled variables.
Rather than introducing a new notation for the variables, 
we choose, wherever convenient, to write, for instance, $x \mapsto x/\delta$. 
This is to be read ``$x$ is replaced by $x/\delta$'', 
so that hereafter the symbol $x$ will denote a scaled variable. 
With this understanding we define
\begin{equation}\label{def:indep-scale}
x\mapsto x/\delta \quad\textrm{and}\quad t\mapsto t/\delta
\end{equation}
and 
\begin{equation}\label{def:uv-scale}
u \mapsto U+\epsilon u\quad\textrm{and}\quad v\mapsto \delta\epsilon v
\end{equation}
so that $U=O(1)$ as $\epsilon\to0$. 
Moreover we write that 
\begin{equation}\label{def:hP-scale}
h=h_0+\epsilon\eta\quad\textrm{and}\quad P=P_{\textrm{atm}}+g(h_0-y)+\epsilon p.
\end{equation}
Physically, $\epsilon\eta$ is the surface displacement from the undisturbed fluid depth $h_0$ 
and $\epsilon p$ is the dynamic pressure, measuring the deviation from the hydrostatic pressure;
see \cite[Section~1.3.2 and Section~3.4.1]{Johnson1997} for further details.

Inserting \eqref{def:indep-scale}--\eqref{def:hP-scale} into \eqref{E:euler1}--\eqref{E:top-bdry2}, we arrive at that:
\begin{equation}\label{E:ww-scale} \left\{
\begin{aligned}
&u_t+Uu_x+U'v+\epsilon(uu_x+vu_y)=-p_x, \\
&\delta^2(v_t+Uv_x+\epsilon(uv_x+vv_y))=-p_y \qquad&&\textrm{in }\{0<y<h_0+\epsilon\eta\}, \\
&u_x+v_y=0, \\
&v=\eta_t+U\eta_x+\epsilon u\eta_x\quad\textrm{and}\quad p=g\eta &&\textrm{on }\{y=h_0+\epsilon\eta\},\\
&v=0 &&\textrm{on }\{y=0\}.
\end{aligned} \right. \end{equation}
Note that $u$, $v$, $p$, $\eta=0$, which corresponds to no deviation from the background current, does satisfy \eqref{E:ww-scale} for arbitrary $h_0>0$ and $U\in C^1([0,h_0])$. 

The governing equations may be nondimensionalized (see, e.g, \cite[Section~1.3.1]{Johnson1997}), 
but we do not pursue that approach here. Rather, we work in the physical variables as much as possible. 

\subsection{Traveling waves and the bifurcation condition}
Waves in the ocean are often observed to be regularly shaped and propagating at constant speed in a fixed direction after they disperse away from the wave generating region. Such wave patterns are representative of \emph{traveling} (or \emph{steady}) \emph{waves}. They are solutions of \eqref{E:euler1}--\eqref{E:top-bdry2} 
for which $u$, $v$ and $P$ are functions of $(x-ct,y)$ 
while  $h$ is a function of $x-ct$ for some $c>0$, called the wave speed. 
It follows that $\psi$ displays the same dependence in $x$ 
and each $\Omega(t)$ is a translation of a time independent domain $\Omega$.

We will further assume that there is no horizontal stagnation:
\be u - c < 0 \qquad \textrm{in } \Omega.  \label{nostagnation} \ee
This amounts to excluding from consideration the so-called ``extreme waves." One important consequence of \eqref{nostagnation} is that each streamline is given by the graph of a single-value function defined on the whole real line. Together with the last equation in \eqref{E:euler1}, \eqref{def:vorticity}, and \eqref{def:psi}, this further implies that $\omega$ and $\psi-cy$ are functionally dependent. One may then write that 
\begin{equation}\label{def:gamma}
-\Delta\psi(x-ct,y)=\omega(x-ct,y)=:\gamma(\psi(x-ct,y)-cy)
\end{equation}
for some function $\gamma$ called the vorticity-stream function.  

We call a traveling wave solution of \eqref{E:euler1}--\eqref{E:top-bdry2} a \emph{Stokes wave} 
provided that it is periodic in the $x$-direction 
and there exists a vertical axis with respect to which $u$ and $h$ are even while $v$ is odd. 
In fact, for an arbitrary function $\gamma$, and under some mild conditions, 
all periodic traveling waves of \eqref{E:euler1}--\eqref{E:top-bdry2} are 
symmetric about their crests; see, e.g., \cite{Hur2007symm} and \cite{CEW2007symm}. 

Studies of such waves were initially made by Stokes \cite{Stokes}, 
who envisioned irrotational periodic waves propagating in an infinitely-deep body of water that is at rest at great depths. 
His deductions were far-reaching, albeit formal, and many were ultimately verified by rigorous analysis.  In fact, the existence and qualitative properties of Stokes waves with zero vorticity are by now quite advanced. This success is strongly indebted to the fact that, for irrotational flows, $\psi$ is a harmonic function, which allows the use of tools such as conformal mappings.  Some excellent surveys are \cite{Toland1996Stokes, OS2001book, Groves2004survey, Strauss2010}.
 
On the other hand, when the flow is rotational,  $\psi$ is no longer harmonic but a solution of \eqref{def:gamma}.  Consequently, one must simultaneously treat the dynamics in the fluid region and the problem of determining the boundary. It was not until recently that 
Constantin and Strauss \cite{CS2004} established the existence theory of large-ampltiude waves for general vorticities.  
Specifically, for each $c>0$ and $k>0$ and for an arbitrary   
$\gamma$ of class $C^{1+\alpha}$, for some $\alpha \in (0,1)$, that satisfies a certain abstract ``bifurcation condition", 
they constructed a global continuum of $2\pi/k$-periodic Stokes waves with the vorticity function $\gamma$.  
 This stimulated a great deal of research activity on water waves with vorticity, leading to a number of important results in a wider hydrodynamical context.  Among them, let us mention only a few that are most relevant to the interests of the present work.  In \cite{Hur2006Stokes1, Hur2011Stokes2}, one of the authors studied waves in deep water; in \cite{Wahlen2009critical, EEW2011critical, CSV2014},  the possibility of critical layers was addressed; in \cite{Varvaruca2009extreme, VW2012extreme}, the formation of a corner was studied in an extreme wave; in  \cite{constantin2011discontinuous}, the result in \cite{CS2004} was generalized to allow for discontinuous vorticity.  

The families of vortical Stokes waves in \cite{CS2004} are constructed via bifurcation theory:  Constantin and Strauss show that a curve of small-ampltiude waves branches from the family of shear flows and continues into the large-ampltiude regime, eventually terminating at (but not including) an extreme wave with a stagnation point.  The bifurcation condition is both necessary and sufficient for the initial branching to occur.  It specifically asks that zero be a (simple) generalized eigenvalue of the associated linearized problem about some shear flow. We note that Constantin and Strauss carry out their analysis in the semi-Lagrangian coordinate system $(x,-\psi)$, which has the effect of transforming  the free boundary problem into one set on a fixed domain.  They also select the wave speed $c$ and the vorticity function $\gamma$ at the outset,  
whereas the shear flow $U$ where bifurcation takes place is determined by $\gamma$ and $c$, and the fluid depth $h_0$ varies along the branch of solutions. Here we derive explicit relations between the pressure and the surface wave in physical variables, and hence it is convenient to consider the situation where the shear flow and the fluid depth are be fixed but the wave speed is chosen to ensure that bifurcation occurs.  In this case, the bifurcation condition becomes the following:  we require that
\begin{equation}\label{E:bifur}\left\{
\begin{split}
&(U-c)(\phi''-k^2\phi)-U''\phi=0\qquad\textrm{for }0<y<h_0,\\
&\phi'(h_0)=\Big(\frac{g}{(U(h_0)-c)^2}+\frac{U'(h_0)}{U(h_0)-c}\Big)\phi(h_0)
\quad\textrm{and}\quad\phi(0)=0
\end{split}\right.
\end{equation}
admits a nontrivial solution for some ${\displaystyle c>\max_{0\leq y\leq h_0}U(y)}$ and $k>0$
and that the solution set is one-dimensional (see \cite{HL2008} for details).
The ordinary differential equation in \eqref{E:bifur} is known as the Rayleigh equation or the inviscid Orr--Sommerfeld equation. 

For general vorticities, \eqref{E:bifur} can only be solved numerically; see, e.g.,\cite{Pete2012dispersion}. Nevertheless, in some cases, explicit solutions can be found and the bifurcation condition agrees with dispersion relations; 
see Example~\ref{eg_irrotational}, Example~\ref{eg_constant vorticity} and Example~\ref{eg_discont vorticity}.
For some range of shear flows, moreover, the bifurcation condition may be verified using ODE theory. If $U\in C^2([0,h_0])$, $U''(h_0)<0$ and $U(h_0)>U(y)$ for $0<y<h_0$, for instance, bifurcation takes place for some $c>\max U$ for all $k>0$; see \cite[Lemma~2.5]{HL2008}.   

We pause to remark that in general the bifurcation condition and the dispersion relation are different: The bifurcation condition is necessary and sufficient to ensure that small-amplitude nontrivial solution exists, whereas the dispersion relation is only a necessary condition ensuring that a plane wave solution exists to the associated linear problem.  At a bifurcation point, one may employ the implicit function theorem (or, more accurately, a Lyapunov--Schmidt reduction procedure), and thus it is possible to find all nontrivial solutions nearby as well as compute them to any order, though perhaps not explicitly. This makes estimating the error quite straightforward. 

In the long wave limit as $k\to0+$, note that the bifurcation condition leads to the Burns condition 
\begin{equation}\label{E:Burns}
\int^{h_0}_0 \frac{dy}{(U(y)-c)^2}=\frac1g;
\end{equation}
see, e.g., \cite{Burns, Johnson1997, HL2008}.


A \emph{solitary water wave} means a traveling wave solution of \eqref{E:euler1}--\eqref{E:top-bdry2}, 
for which $h(x-ct)$ tends to a constant and $v(x-ct,y)\to 0$ uniformly for $y$ as $x-ct\to\pm\infty$. 
It may be viewed formally as the limit of Stokes waves as the period tends to infinity. 
Indeed, small-amplitude solitary water waves, when they exist, emanate  
near the critical speed determined upon solving \eqref{E:Burns}. 
Solitary water waves are yet more technically challenging to study because they spread over an unbounded domain. 
Classical bifurcation theory is less useful as a consequence.  
(For periodic waves, one can imagine $x$ as lying on a circle, which compactifies the problem).  
In the case of zero vorticity, the rigorous existence theory goes back to 
the constructions in \cite{FH1954, Beale1977solitary} of small-ampltiude waves 
and it includes singular bifurcation result in \cite{AT1981solitary} of large-ampltiude waves. 
Recently, these results have been extended to include an arbitrary distribution of vorticity 
in \cite{Hur2008solitary, GW2008solitary,Wheeler2013solitary,wheeler2015pressure}. 

\subsection{The pressure transfer function}\label{sec:ww-transfer}
In the linear approximation, suitable for small-ampltiude (but not necessarily long wavelength) waves, 
$\epsilon\ll1$ while $\delta = O(1)$. 
Under this assumption, \eqref{E:ww-scale} becomes to leading order:
\begin{subequations}\label{E:ww-linear} 
\begin{equation}\label{E:euler-linear}\left\{
\begin{split}
&u_t+Uu_x+U'v=-p_x,\qquad&&\\
&v_t+Uv_x=-p_y \qquad&&\textrm{in }\{ 0<y<h_0 \}, \\
&u_x+v_y=0,
\end{split}\right.
\end{equation}
and
\begin{align}
&v=\eta_t+U\eta_x\quad\textrm{and}\quad p=g\eta &&\text{on }\{y=h_0\},\label{E:top-linear} \\
&v=0 &&\text{on }\{y=0\}.\label{E:bed-linear}
\end{align}\end{subequations}
Here, for simplicity, we have taken $\delta = 1$.  

Note that small-amplitude Stokes waves with vorticity constructed in \cite{CS2004} solve \eqref{E:ww-linear} with errors of $O(\epsilon^2)$ (in $C^{2+\alpha}$ for the velocity and the pressure, and in $C^{1+\alpha}$ for the free surface). 
Moreover 
the wave speed is approximated by that in \eqref{E:bifur} with errors of $O(\epsilon^3)$.

Seeking $2\pi/k$-periodic solutions of \eqref{E:ww-linear} traveling at the speed $c>0$ and with the wave number $k>0$, it is reasonable to assume that 
\begin{equation}\label{E:Stokes-ansatz}
\eta(x;t)=\cos(k(x-ct))
\end{equation}
to leading order. 
One can then solve \eqref{E:ww-linear} explicitly to find that
\begin{equation}\label{E:linear-soln}
\begin{split}
u(x,y;t)&=\frac1k\phi'(y)\cos(k(x-ct)), \\
v(x,y;t)&=\phi(y)\sin(k(x-ct)),\\
p(x,y;t)&= \frac1k((c-U(y)\phi'(y)+U'(y)\phi(y))\cos(k(x-ct)),
\end{split}\end{equation}
where $\phi$ solves  
\begin{equation}\label{E:Rayleigh} \left\{
\begin{split}
&(U-c)(\phi''-k^2\phi)-U''\phi=0\qquad \textrm{for }0<y<h_0,\\
&\phi(h_0)=k(c-U(h_0))\quad\textrm{and}\quad\phi(0)=0.
\end{split} \right. 
\end{equation}
It is easy to see that the bifurcation condition \eqref{E:bifur} can be deduced from \eqref{E:Rayleigh} 
by employing the first equation in \eqref{E:top-linear}. 
In particular, the solutions of \eqref{E:Rayleigh} is unique.
Indeed, a non-trivial solution of \eqref{E:bifur} is unique up to multiplication by a constant. 
The normalization in \eqref{E:Stokes-ansatz} would then single out the unique solution of \eqref{E:bifur}, and hence \eqref{E:Rayleigh}.

The last equation in \eqref{E:linear-soln} allows us to define the linear {\em pressure transfer function}:
\begin{equation}\label{def:transfer}
T(y)=\frac1k((c-U(y))\phi'(y)+U'(y)\phi(y)),
\end{equation}
which relates the dynamic pressure and the surface displacement from the undisturbed fluid depth via the identity
\[
p(x,y;t)=T(y)\eta(x;t).
\]
In particular, the displacement of the free surface from its rest height 
and the trace of the dynamic pressure at the bed are related by the formula
\begin{equation}\label{E:bottom-to-top}
\eta(x;t)=\frac{1}{T(0)}p(x,0;t)=k\frac{1}{c-U(0)}\frac{1}{\phi'(0)}p(x,0;t).
\end{equation}
We pause to remark that $\phi'(0)\neq0$, and hence the formula makes sense.
Otherwise, the unique solution of the initial value problem for the Rayleigh equation 
and the second boundary condition in \eqref{E:Rayleigh} would force $\phi \equiv 0$, 
which would then contradict the first boundary condition in \eqref{E:Rayleigh}.

We emphasize that \eqref{def:transfer} is a linear relation.
Indeed, the lack of explicit solutions prevents a formula of the kind for large-ampltiude waves. 
A nonlinear relation between the trace of the pressure and the surface wave is, of course, far more involved.
Nevertheless, it was achieved recently for Stokes and solitary water waves in the case of zero vorticity \cite{Constantin2012,clamond2013recovery,clamond2013new} or constant vorticity \cite{OVDH2012}.
These results make strong use of the fact that one can explicitly solve the Laplace equation in a strip using separation of variables,
and thus they do not extend to rotational flows.
One may instead use higher order perturbations of the governing equations 
to obtain pressure transfer functions with improved accuracy; see, e.g., \cite{HL2015pressure}.

In general, \eqref{E:Rayleigh} must be investigated numerically.  However, for certain special cases the calculations are considerably simplified and analytic methods are possible.  We present a few of these below.  

\begin{example}[Zero vorticity]\label{eg_irrotational}\rm
Suppose that $U\equiv0$, namely that there is zero vorticity.  Then,
a straightforward calculation reveals that a nontrivial solution to \eqref{E:bifur} exists provided that
\begin{equation}\label{E:c0}
c^2=\frac{g\tanh(kh_0)}{k}.
\end{equation}
This is the well-known dispersion relation for gravity water waves in finite-depth.
Moreover the solution to \eqref{E:Rayleigh} is 
\[
\phi(y)=kc\frac{\sinh(ky)}{\sinh(kh_0)}. 
\]
Inserting this and \eqref{E:c0} into \eqref{def:transfer}, we find the pressure transfer function:
\begin{equation}\label{E:0transfer}
T(y)=g\frac{\cosh(ky)}{\cosh(kh_0)},
\end{equation}
which agrees with the result in \cite{ES2008}, for instance.
In this sense, \eqref{def:transfer} generalizes the well-known pressure transfer function 
to allow for the effects of vorticity.
\end{example}

\begin{example}[Constant vorticity]\label{eg_constant vorticity}\rm
Let $U(y)=\gamma(y-h_0)$, $0\leq y\leq h_0$, for some constant $\gamma$. 
(More generally, one may take $U(y)=\gamma(y-h_0)+U_0$ for some constant $U_0$, 
but $U_0$ can be absorbed to the wave speed.) This models the constant vorticity $-\gamma$. 
A straightforward calculation reveals that a nontrivial solution to \eqref{E:bifur} exists, provided that 
\begin{equation}\label{E:cgamma}
c = c(\gamma) =-\frac{\gamma\tanh(kh_0)}{2k}+\sqrt{\frac{\gamma^2\tanh^2(kh_0)}{4k^2}+\frac{g\tanh(kh_0)}{k}};
\end{equation}
see also \cite{CS2004,vasan2014pressure}. 
The other solution with the $-$ sign violates $c>0$, and hence we discard it. We note in passing that 
\[
c(-|\gamma|) > c(0) > c(|\gamma|),
\]
meaning that favorable currents corresponding to $\gamma<0$ enhance the wave speed, 
while adverse currents ($\gamma<0$) reduce it. 
Moreover one may explicitly solve \eqref{E:Rayleigh} to find that
\[
\phi(y)=kc\frac{\sinh(ky)}{\sinh(kh_0)}.
\] 
This leads to the pressure transfer function: 
\begin{equation}\label{E:ctransfer}
T(y)=\frac{c}{\sinh(kh_0)}((c-\gamma(y-h_0))k\cosh(ky)+\gamma\sinh(ky)),
\end{equation}
where $c$ is in \eqref{E:cgamma}. When $\gamma=0$, clearly, \eqref{E:ctransfer} becomes \eqref{E:0transfer}.
\end{example}

\begin{remark*}\rm
It is readily seen from \eqref{E:cgamma} that a necessary and sufficient condition for a critical layer to exist in a laminar flow is that
\begin{equation}\label{stag point}
\gamma < 0\quad \text{and } \quad \gamma^2 > {g\tanh(kh_0) \over kh^2_0 - h_0\tanh(kh_0)}.
\end{equation}
In the small-ampltiude regime, therefore, a nontrivial Stokes wave with a constant vorticity will contain a stagnation point if and only if \eqref{stag point} is satisfied by the background shear flow; see also \cite{vasan2014pressure}. 
\end{remark*}

\begin{remark*}\rm
One may deduce from Example~\ref{eg_constant vorticity} the effects of vorticity on the pressure distribution. In the case of zero vorticity, it follows from \eqref{E:0transfer} that the dynamic pressure to the leading order is minimized at the bottom and the pressure disturbance attenuates in depth. Note from \eqref{E:ctransfer} that 
\[
T'(y) = \frac{k^2c\sinh(ky)}{\sinh(kh_0)}(c-\gamma(y-h_0)) .
\]
Hence if $\gamma > 0$, corresponding to a negative constant vorticity, then
\[
T(y), T'(y) > 0\quad \text{ for } 0 < y \leq h_0,
\]
leading to a similar behavior of the pressure as in the irrotational case. However, if $\gamma$ is sufficiently negative satisfying \eqref{stag point} then $T'(y)$ may become negative, and the pressure disturbance becomes non-monotone in depth; see, e.g., Figure \ref{const vorticity}.   This agrees with numerical observations made by Ko and Strauss in \cite{ko2008large,ko2008vorticity}.  
\begin{figure}[tb]
  \includegraphics[scale=0.5]{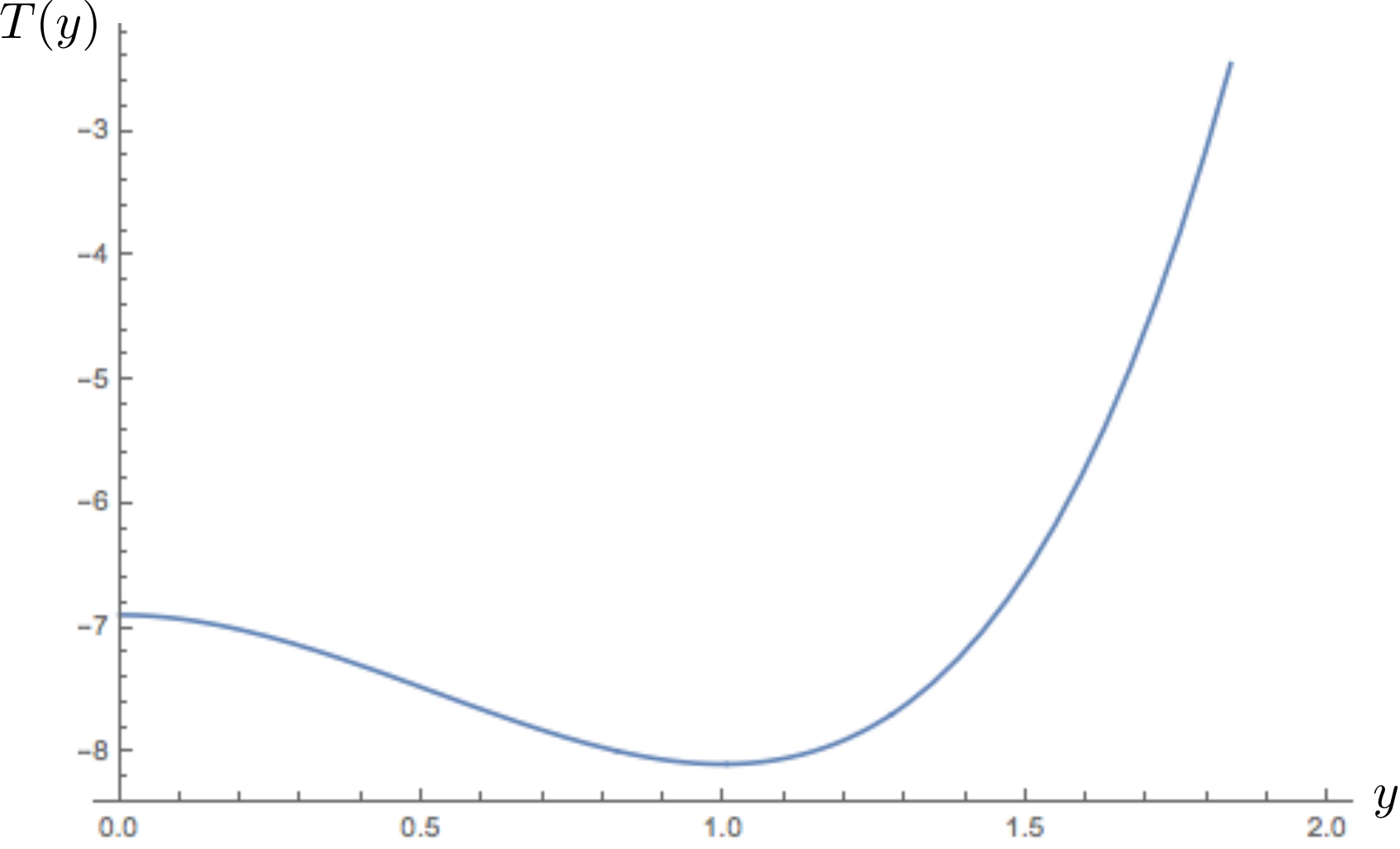}
  \caption{Pressure transfer function $T(y)$ for $\gamma = -5, k = 1, h_0 = 2$.}
  \label{const vorticity}
\end{figure}
\end{remark*}

\begin{example}[Piecewise constant vorticity]\label{eg_discont vorticity}\rm
Let
\[
U(y)=\begin{cases}
\gamma_-(y-h_0) \quad &\textrm{for }0 \leq y \leq h_1,\\
\gamma_-(h_1-h_0)+\gamma_+(y-h_1)& \textrm{for }h_1< y\leq h_0.
\end{cases}
\]
Therefore
\be \label{piecewise constant vorticity} 
U'(y)=\begin{cases}
\gamma_- \quad & \textrm{for }0 \leq y \leq h_1,\\
\gamma_+ & \textrm{for }h_1< y \leq h_0
\end{cases} \quad\textrm{and}\quad U''(y)= (\gamma_+ - \gamma_-) \delta(y-h_1),
\ee
where $\delta$ denotes the Dirac measure with unit mass. 
In other words, the vorticity is piecewise constant. 
The second equation in \eqref{piecewise constant vorticity} means that 
the Rayleigh equation must be understood in the sense of distributions. 

The existence of Stokes waves with discontinuous vorticity was established 
in \cite{constantin2011discontinuous}, subject to a bifurcation condition.  Analogous to the smooth vorticity case, it asks that a certain Sturm--Liouville problem posed in the semi-Lagrangian variables admits a nontrivial distributional solution. 
When the vorticity is piecewise smooth like \eqref{piecewise constant vorticity}, 
it is straightforward to translate this condition back to physical variables. 
Ultimately, one finds that bifurcation occurs if \eqref{E:bifur} has a nontrivial weak solution
in $C^2([0,h_0] \setminus \{h_1\}) \cap W^{1,\infty}([0,h_0])$.   

Let $\phi$ be the solution of \eqref{E:Rayleigh} and 
let $\phi_+$ and $\phi_-$ be its restriction to $(h_1, h_0]$ and $[0,h_1)$, respectively. 
The Rayleigh equation \eqref{E:Rayleigh} with discontinuous coefficients determines $\phi_\pm$ up to two constants.
The continuity of $\phi$ on $[0,h_0]$ implies that
\[ \phi_+(h_1) = \phi_-(h_1),\]
while the distributional identity 
\[ \phi''= \frac{1}{U-c}\delta(\cdot-h_1)\phi\]
gives that
\[  \phi'_+(h_1) - \phi'_-(h_1) = {\gamma_+ - \gamma_- \over c - U_1} \phi(h_1),\]
where
\[
U_1= U(h_1) =\gamma_- (h_1-h_0).
\]
We therefore infer that 
\[
\phi(y)=\begin{cases}
A_-\sinh(ky) \qquad &\textrm{for }0 \leq y \leq h_1, \\ 
A_+\sinh(k(y - h_0)) +B\cosh(k(y - h_0))&\textrm{for }h_1\leq y \leq h_0,
\end{cases}\]
where $B=k(c-U(h_0))=k(c-(\gamma_--\gamma_+)(h_1-h_0))$ and 
\begin{align*}
A_-&={\displaystyle\sinh(k(h_1-h_0))-\cosh(k(h_1-h_0))\coth(k(h_1-h_0))\over 
\displaystyle\sinh(kh_1)\left({\gamma_+-\gamma_-\over k(c - U_1)}+\coth(kh_1)-\coth(k(h_1-h_0))\right)}B,\\
A_+ &={\displaystyle \left({\gamma_+-\gamma_- \over k(c - U_1)}+\coth(kh_1)\right)\coth(k(h_1-h_0))-1 
\over\displaystyle\coth(k(h_1-h_0))-\left({\gamma_+-\gamma_-\over k(c-U_1)}+\coth(kh_1)\right)}B.
\end{align*}
Inserting this into \eqref{E:linear-soln}, incidentally, we obtain the implicit dispersion relation:
\begin{equation}\label{c discont}
(c - U(h_0))(A_++\gamma_+) = g.
\end{equation}

The pressure transfer function can likewise be computed explicitly 
using \eqref{def:transfer}, $\phi$ and the dispersion relation. 
The expression is quite complicated, however, and we do not record it here.
Instead, we provide the transfer function from the pressure at the bed to the surface wave: 
\[
T(0) =\frac1k(c-U(0))\phi'(0) = A_-(c+\gamma_-h_0).
\]
\end{example}

\subsection{The hydrostatic approximation}\label{sec:ww-hydrostatic}
In a weakly nonlinear approximation, suitable for small-amplitude solitary water waves, $\delta^2=\epsilon\ll1$. 
With this scaling regime, \eqref{E:ww-scale} becomes to leading order
\begin{equation}\label{E:ww-wnonlinear}
\left\{
\begin{aligned}
&u_t+Uu_x+U'v=-p_x, \\
&p_y=0 \qquad&&\text{in }\{0<y<h_0\}, \\
&u_x+v_y=0, \\
&y=\eta_t+U\eta_x\quad\text{and}\quad p=g\eta \qquad&&\text{on }\{y=h_0\},\\
&v=0 &&\text{on }\{y=0\}.
\end{aligned} \right.
\end{equation}

Small-amplitude solitary water waves with an arbitrary vorticity found in \cite{Hur2008solitary}, 
for instance, are exact solutions of \eqref{E:euler1}--\eqref{E:top-bdry2}, and hence \eqref{E:ww-scale}.
As $\delta^2=\epsilon\to0$, furthermore, it follows from the existence proof that 
they are approximated to leading order by the famous Korteweg-de Vries (KdV) soliton, 
and the wave speed is given to leading order by  \eqref{E:Burns}.
In other words, small-amplitude solitary water waves with vorticity will solve \eqref{E:ww-wnonlinear} 
with errors of $O(\epsilon^2)$ in the space of real analytic functions.  This regularity assumption is not as stringent as it first appears.   
For an arbitrary vorticity function in the $C^{\alpha}$ class, for some $\alpha \in (0,1)$, 
 one of the authors proved in \cite{Hur2012regularity} that 
any solitary water wave with velocity in  $C^{1+\alpha}$ is in fact real analytic. 

Observing that the second equation in \eqref{E:ww-wnonlinear} implies that 
the dynamic pressure is independent of the depth at leading order, 
we arrive at the hydrostatic approximation: 
\begin{equation}\label{def:ww-static}
\eta(x;t)=\frac1gp(x,0;t).
\end{equation}
It is straightforward to verify that exact relations in \cite{OVDH2012} and \cite{Constantin2012} 
between the trace of the pressure and the surface of a solitary water wave agree with \eqref{def:ww-static} to leading order. 
Note that \eqref{def:ww-static} is independent of the underlying shear flow. 
Higher order approximations, however, incorporate the effects of vorticity;
see \cite{HL2015pressure} for details.

\begin{remark*}\rm  
The study of pressure transfer functions is often motivated by the desire to track tsunamis. 
When working with long wave approximations, however, some caution is warranted.  
As is pointed out in \cite{constantin2006tsunami,constantin2008propagation,constantin2008nondimen},
for instance, the time scale on which the KdV dynamics manifests is 
considerably shorter than the lifespan of a typical tsunami. 
This unfortunately limits the usefulness of \eqref{def:ww-static} as a means of reconstructing tsunamis,
though it may be applicable to other phenomena such as slowly moving waves.  
\end{remark*}

\section{Surface waves with continuous stratification}\label{sec:continuous}

The presence of salinity and temperature gradients may lead to a heterogeneous density distribution that strongly influences the flow.  In this section, we endeavor to study the pressure transfer function in stratified fluids.  
We continue to assume that the fluid is inviscid and incompressible,
and we use the notation in Section~\ref{sec:vorticity}.
The main difference is that the density of the fluid $\rho$ is no longer a constant in \eqref{E:euler1}
and instead it is subject to the equation of mass conservation:
\begin{equation}\label{E:density}
\rho_t + u\rho_x + v\rho_y = 0 \quad \textrm{in }\Omega(t).
\end{equation}
In order for the flow to be physically realistic, 
we assume that $\rho$ is a non-negative function and non-increasing in the $y$-variable. 
The boundary conditions are exactly as stated in Section~\ref{sec:vorticity}. 
Note that \eqref{E:density} implies that the density is constant on the free surface and the bed for traveling waves. 

The existence of periodic traveling waves of \eqref{E:euler1}--\eqref{E:top-bdry2}, \eqref{E:density} 
has been the subject of much research, dating back to the early work of Dubreil-Jacotin \cite{dubreil1937theoremes}.  Some notable results include \cite{amick1984semilinear,terkrikorov1963theorie,turner1981internal,turner1984variational,lankers1997fast}. 
The result most relevant to our purposes is \cite{walsh2009stratified} on the existence large-amplitude 
 periodic traveling waves with general density stratifications and general background currents. 
Specifically, for each $c>0$ and $k>0$ and for an arbitrary function of class $C^{1+\alpha}$, for some $\alpha\in(0,1)$, relating the stream function and density stratification and subject to a bifurcation condition, 
there exists a global continuum of $2\pi/k$-periodic traveling waves of \eqref{E:euler1}--\eqref{E:top-bdry2}, \eqref{E:density}.

Less is known about solitary waves, because the unboundedness of the domain leads to serious issues with compactness. 
All existing literature on the existence of solitary waves in stratified fluids pertains exclusively to the case of zero background current, and no exact existence theory is available for general density stratifications and general shear flows in the far field. This is the subject of a forthcoming paper by two of the authors.  
For now we merely focus on periodic traveling waves.

Note that
\begin{equation}\label{basic state}
h\equiv h_0, \quad (u, v) = (U(y), 0), \quad \rho = R(y)
\quad \textrm{and} \quad  P = P_{\textrm{atm}} + \int^{h_0}_y g{R}(z)~dz
\end{equation}
form a solution of \eqref{E:euler1}-\eqref{E:top-bdry2}, \eqref{E:density}
for arbitrary $h_0>0$, $U\in C^1$, and $R(y)\geq0$, $R'(y)\leq 0$ for all $0<y<h_0$. 
We are interested in waves propagating over the shear flow and the density stratification of the form. 
In the remainder of the section, therefore, $h_0$, $U$ and $R$ are held fixed. 

To proceed, we repeat the procedure in Section~\ref{ss:scaling} and 
introduce the long wavelength parameter $\delta$ and the small-ampltiude parameter $\epsilon$. 
We define scaled variables $x$, $t$ and $u$, $v$ by \eqref{def:indep-scale}--\eqref{def:uv-scale}, and we write that
\begin{equation}\label{Stratified scaling} 
h=h_0 + \epsilon \eta, \qquad  P=P_{\textrm{atm}} + \int^{h_0}_y g R(z)~dz + \epsilon p
\quad\textrm{and}\quad \rho=R + \epsilon \rho.
\end{equation} 
Physically, $\epsilon\rho$ measures the deviation from the background density. 

Inserting \eqref{def:indep-scale}--\eqref{def:uv-scale} and \eqref{Stratified scaling} into \eqref{E:euler1}--\eqref{E:top-bdry2}, \eqref{E:density},
we arrive, in the linear approximation as $\epsilon\to0$ and $\delta=1$, at that:
\begin{equation}\label{Stratified linear} \left\{
\begin{aligned}
&\rho_t + U\rho_x + v R' = 0,\\
&\displaystyle R \left( u_t + Uu_x + U' v \right) = -{p_x}\qquad&& \text{in } \{ 0< y < h_0 \},\\
&\displaystyle R (v_t + U v_x) = -{ p_y -  \rho },\\
&u_x + v_y = 0,\\
&v = \eta_t + U\eta_x \quad \text{and}\quad p = g R \eta\qquad&& \text{on }\{y = h_0\}, \\
&v = 0\qquad&& \text{on }\{y = 0\}.
\end{aligned} \right.
\end{equation}
Periodic traveling waves found in \cite{walsh2009stratified}, for instance, are exact solutions of \eqref{E:euler1}--\eqref{E:top-bdry2}, \eqref{E:density}, and hence approximate solutions of \eqref{Stratified linear}. 
From the bifurcation theory analysis, one deduces that, if $\rho \in C^{1+\alpha}$, for some $\alpha \in (0,1)$, then \eqref{Stratified linear} holds up to an error term of order $O(\epsilon^2)$ in $C^{2+\alpha}$ for the velocity, $C^{3+\alpha}$ for the pressure, and $C^{1+\alpha}$ for the density.   

With that in mind, we again assume that $\eta$ has the ansatz
\[
\eta(x;t)=\cos(k(x-ct)),
\]
which allows us to solve \eqref{Stratified linear} explicitly.  We find that 
\be \label{linear stratified solution} 
\begin{split}
u(x,y; t) & = \frac{1}{k} \phi^\prime(y) \cos{(k(x-ct))}, \\
v(x,y; t) & = \phi(y) \sin{(k(x-ct))}, \\
p(x,y; t) & = \frac1k R(y) \left((c-U(y))\phi'(y)+U'(y)\phi(y) \right)\cos(k(x-ct)),\\
\rho(x,y; t) & = \frac{1}{c-U} R^\prime(y) \phi(y) \sin{(k(x-ct))},
\end{split}
\ee where $\phi$ solves
\begin{equation}\label{Stratified Rayleigh1} \left\{
\begin{aligned}
&k^2(c-U)^2 R \phi=(c - U)(R((c - U) \phi' + U'\phi))' - R' \phi \quad \text{for }0< y < h_0,\\
&\phi(h_0) = k( c - U(h_0))\quad\textrm{and}\quad \phi(0) = 0.
\end{aligned} \right.
\end{equation}
When $R\equiv1$, clearly, this reduces to \eqref{E:Rayleigh}. 

Under the assumption ${\displaystyle c>\max_{0\leq y\leq h_0}U(y)}$, namely no critical layers,
\eqref{Stratified Rayleigh1} may be transformed into a more amenable form by making change of unknowns.  Let 
\begin{equation*}\label{change of unknown}
\varphi= {\phi \over c - U},
\end{equation*}
and \eqref{Stratified Rayleigh1} becomes
\begin{equation}\label{Stratified Rayleigh2} \left\{
\begin{aligned}
&(R (c - U)^2 \varphi')'=k^2(c - U)^2 R + R')\varphi \quad \text{for }0< y < h_0,\\
&\varphi(h_0) = k \quad\textrm{and}\quad \varphi(0) = 0.
\end{aligned} \right.
\end{equation}
Ultimately, we define the linear pressure transfer function:
\begin{equation}\label{transfer stratified}
T(y)=\frac1k R(y)((c-U(y))\phi'(y)+U'(y)\phi(y)) = \frac1k R(y)(c-U(y))^2\varphi'(y),
\end{equation}
which relates the dynamic pressure and the surface displacement from the undisturbed fluid depth 
via the identity
\begin{equation}\label{pressure stratified}
p(x,y;t)=T(y)\eta(x;t).
\end{equation}

\begin{example}[Exponential stratification]\rm
It is in general impossible to solve either \eqref{Stratified Rayleigh1} or \eqref{Stratified Rayleigh2} explicitly.  An example for which one may find a solution in closed form is 
when the density increases exponentially with depth, i.e.
\[
R(y) = e^{-2\beta y} \quad \text{for some } \beta>0.
\]
Then,  \eqref{Stratified Rayleigh2} becomes 
\[
((c-U)^2\varphi')' - 2\beta(c-U)^2\varphi' = (k^2(c-U)^2 - 2\beta)\varphi.
\]

Consider the simplest case $U\equiv0$. 
(More generally one may take $U$ to be a constant, but this constant can just be absorbed into the wave speed.)
The above equation then reduces to the second-order constant-coefficient ODE: 
\begin{equation}\label{ode}
c^2\varphi'' - 2\beta c^2 \varphi' + (2\beta - k^2c_1^2) \varphi = 0, \qquad 
\varphi(0) = 0\quad\textrm{and}\quad \varphi(h_0) = k,
\end{equation}
and the sign of the discriminant of the characteristic equation $4c^2((\beta^2 + k^2)c^2 - 2\beta)$ 
determines the type of the solution. 
Furthermore, applying \eqref{transfer stratified} and \eqref{pressure stratified} 
to the linearized kinematic equation in \eqref{Stratified linear}, we determine that the dispersion relation is
\[
c^2= {gk \over \varphi'(h_0)}.
\]
Note that this is in fact well-defined. Indeed, it is easily seen that the first boundary condition in \eqref{Stratified Rayleigh2} can  be written as
\[
\varphi'(h_0) = {g\over c^2} \varphi(h_0).
\]
Since $\varphi(h_0) = k$, it follows that $\varphi'(h_0) \neq 0$.
%
We then infer from \eqref{transfer stratified} the linear pressure transfer function:
\[
T(y) = {1\over k} c^2e^{-\beta y} \varphi'(y) = ge^{-\beta y} {\varphi'(y) \over \varphi'(h_0)}.
\]
\end{example}

\section{Internal waves and wind-driven waves}\label{sec:internal}
We now turn our attention to traveling waves propagating 
at the interface between two immiscible fluids in a channel.  
Suppose that the fluid domain is partitioned into two layers and 
confined to a strip which is bounded from above by a flat rigid lid at $\{y=H\}$ for some $H \in (h_0,\infty]$
and from below by a flat rigid bed at $\{y = 0\}$. 
We assume that the two layers share a common boundary 
that is given as the graph of a smooth function $h = h(x;t)$ satisfying $0<h<H$. 
Let $\Omega(t) = \Omega_-(t) \cup \Omega_+(t)$, where
\begin{align*}  \Omega_+(t) & := \{  (x,y) \in \mathbb{R}^2 :  h(x;t) <  y < H \}, \\
\Omega_-(t) & := \{ (x, y) \in \mathbb{R}^2 :  0 < y < h(x;t) \},
\end{align*}
denote the upper and lower fluid layers, respectively, and let 
\[  S(t) = \partial\Omega_+(t)\cap \partial\Omega_-(t) = \{ y = h(x;t) \}\]
be the internal interface separating them.
For the velocity, pressure, and density, we will use the notation in Section~\ref{sec:vorticity} and Section~\ref{sec:continuous}.

The velocity field and pressure satisfy \eqref{E:euler1} in $\Omega(t)$.
We assume for simplicity that the density is layer-wise constant:
\[ \rho(x,y;t) = \rho_+ \chi_{\Omega_+(t)}(x,y) + \rho_- \chi_{\Omega_-(t)}(x,y)\]
for some fixed $\rho_\pm > 0$, 
where $\chi_{\Omega_\pm(t)}$ denotes the characteristic function for $\Omega_\pm(t)$.  
Note that \eqref{E:density} holds trivially in $\Omega(t)$.  
Here and in the sequel, $\pm$ denotes the restriction to $\Omega_\pm(t)$ of a function defined on $\Omega(t)$.  

The motion of the free interface is governed by the kinematic and dynamic conditions:
\be
v_\pm = h_t+u_\pm h_x \quad \text{and} \quad P_+ - P_- = \sigma \frac{h_{xx}}{\left(1+h_x^2 \right)^{3/2}}  \quad \textrm{on } S(t), \label{E:top-bdry4} \ee
where $\sigma \geq 0$ is the coefficient of surface tension.  
The first equation is equivalent to requiring that the interface is a material line, or 
that the normal component of the velocity field is continuous on $S(t)$.  
The second comes from the Young--Laplace law, stating that on each point of the free interface, 
the pressure experiences a jump that is proportional to the (signed) curvature of $S(t)$ there.  
When $\sigma = 0$, this of course means that the pressure is continuous in $\overline{\Omega(t)}$.  
For now we include surface tension because the problem is ill-posed otherwise; 
see, e.g., \cite{shatah2011interface}.  
We impose an impermeability condition on the bed and the lid:
\be v = 0 \qquad \text{on }\{y = 0 \} \cup \{ y = H \}. \label{E:b-bdry4} \ee  

Many stratified fluids are observed to have layers of constant density 
separated by thin pycnoclines, wherein the density undergoes a rapid change. 
It is convenient for computations to simply neglect the thickness of the pycnocline, 
which at least formally leads to a two- or multiple-fluid approximation.  
This can in fact be made rigorous in some situations; 
see, e.g., \cite{chen2015dependence} and \cite{james2001internal}.  

Another circumstance one encounters the two-fluid setup is in the theory of wind generation of water waves,
for which $\Omega_+(t)$ represents the air region and $\Omega_-(t)$ is the water;
see, e.g., \cite{janssen2004interaction} for more discussion.
Note that, in contrast to the previous sections, this model takes into account the dynamics in the atmosphere as well as the water.  
The upper boundary is a mathematical contrivance that is reasonable 
because the flow in the atmosphere is expected to have negligible influence on the motion of the free surface.  

The existence of periodic traveling waves and solitary waves to  
\eqref{E:euler1}, \eqref{E:density}, \eqref{E:top-bdry4}, \eqref{E:b-bdry4} has been established
in a number of works, although a majority is restricted to 
the absence of the background flow and the effects of surface tension; 
see, e.g., \cite{turner1981internal,turner1984variational,amick1986global}.  
For a nontrivial background current, but neglecting capillary effects, 
one of the authors in \cite{walsh2013wind} constructed a continuum of small-amplitude solutions.  

Note that 
\begin{equation}\label{E:background}
h\equiv h_0, \qquad (u,v)=(U_\pm(y),0) \quad\textrm{and}\quad P=g\rho(h_0-y)
\end{equation}
form a solution of \eqref{E:euler1}, \eqref{E:density}, \eqref{E:top-bdry4}, \eqref{E:b-bdry4} 
for arbitrary $0<h_0<H$ and $U_\pm$ in the $C^1$ class.  These are the family of shear flows that we take as the model for the background current.  

Repeating the procedure in Section~\ref{ss:scaling}, we introduce \eqref{def:parameters}
and we scale the independent and dependent variables using \eqref{def:indep-scale}--\eqref{def:hP-scale}, but 
\be P =g\rho(h_0 - y) + \epsilon\rho p \label{def dynamic pressure} \ee
for the pressure. Substituting these, \eqref{E:euler1} becomes
\begin{equation}\label{E:euler4-scaled} \left\{
\begin{aligned}
&u_t+Uu_x+U'v+\epsilon(uu_x+vu_y)=-p_x, \\
&\delta^2(v_t+Uv_x+\epsilon(uv_x+vv_y))=-p_y \qquad&&\textrm{in }\Omega(t), \\
&u_x+v_y=0, \\
\end{aligned} \right. \end{equation}
while \eqref{E:top-bdry4} turns in to 
\be \label{scaled kinematic cond} 
v_\pm = \eta_t + U_\pm \eta_x + \epsilon u_\pm \eta_x  \qquad \textrm{on } S(t) \ee 
and 
\be \label{scaled jump cond} 
p_- - \frac{\rho_+}{\rho_-} p_+  = \left(1 - \frac{\rho_+}{\rho_-} \right) g \eta - \frac{\sigma \delta^2 }{\rho_-} \frac{\eta_{xx}}{\left( 1+ \epsilon^2 \delta^2 \eta_x^2 \right)^{3/2}} \qquad \textrm{on } S(t). \ee
The bottom boundary condition is exactly the same as in \eqref{E:b-bdry4}.

\subsection{The pressure transfer function}

In the linear approximation, suitable for small-amplitude periodic traveling waves, 
we repeat the procedure in Section~\ref{sec:ww-transfer}:  formally take $\delta=1$ and $\epsilon=0$ in  
\eqref{E:euler4-scaled}-\eqref{scaled jump cond} and assume that 
\[
\eta(x;t)=\cos(k(x-ct)).
\]
A straightforward calculation then reveals that
\begin{align} \label{two fluid stokes}
u_\pm(x,y;t) &=\frac1k (c-U_\pm(h_0)) \phi_\pm'(y)\cos(k(x-ct)),\notag \\
v_\pm (x,y;t) &= (c-U_\pm(h_0)) \phi(y) \sin(k(x-ct)), \\
p_\pm(x,y;t) &= T_\pm(y) \cos{(k(x-ct))},\notag
\end{align}
where 
\be \label{two fluid stokes transfer} 
T_\pm(y)=  (U_\pm(y) - c)(U_\pm(h_0)-c) \phi_\pm'(y) - U_\pm'(y)(U_\pm(h_0)-c) \phi_\pm(y) \ee
and $\phi_\pm$ solves
 \be \label{rayleigh} \begin{cases}
- \phi_\pm'' + \left(\dfrac {U_\pm''}{U_\pm-c} + k^2\right) \phi =0\qquad \text{for }y\in(0,h_0)\cup(h_0, H), \\
 \phi_-(0) = \phi_+(H) =0\quad\textrm{and}\quad \phi_\pm(h_0)=k. \end{cases}
\ee
The pressure transfer function is defined by \eqref{two fluid stokes transfer}.

We mention that evaluating \eqref{scaled jump cond} with $\epsilon = 0$ and $\delta = 1$, gives
\[ p_- - \frac{\rho_+}{\rho_-} p_+ = \left( \left(1- \frac{\rho_+}{\rho_-} \right) g + k^2 \frac{\sigma}{R_-} \right) \cos{(k(x-ct))} \qquad \textrm{on }\{ y = h_0 \}.  \]
Combining this with the last equation in \eqref{two fluid stokes} and \eqref{two fluid stokes transfer}, 
we arrive at the dispersion relation:  
\be \label{two fluid stokes dispersion}
\begin{split} 
(R_+ - R_-) g - k^2 \sigma &=  R_+ T_+(h_0) - R_- T_-(h_0) \\
&=   R_+(U_+(h_0) - c)^2 \phi_+'(h_0) - kR_+U_+'(h_0)(U_+(h_0)-c) \\
& \quad -  R_- (U_-(h_0) - c)^2 \phi_-'(h_0) + k R_- U_-'(h_0)(U_-(h_0)-c). 
\end{split}
 \ee
Therefore the effects of surface tension enter the pressure transfer function through the wave speed,
which is determined by the dispersion relation. 

For a general background current, $\phi_\pm$ will not be found in closed form.  
Nonetheless, it is relatively simple to show that in the case of $U \equiv 0$ 
this reduces to the well-known dispersion relation for two-fluid irrotational flow (see, e.g., \cite{buhler2015wind}).  
It is worth mentioning that formally as $R_+/R_- \to 0$, 
we recover \eqref{E:c0} if $U \equiv 0$, and \eqref{E:cgamma} if $U$ is linear.   

\subsection{The hydrostatic approximation}
In a small-amplitude solitary wave regime, where $\delta^2=\epsilon\to0$ , 
the second equation in \eqref{E:euler4-scaled} reduces to
\be p_y = 0 \qquad \textrm{in } \Omega,\label{vertical momentum 0} \ee
where recall that $\Omega = \mathbb{R} \times (0,H) \setminus \{ y = h_0\}$,
and \eqref{scaled jump cond} becomes 
\be p_{-} - \frac{\rho_+}{\rho_-} p_{+} = \left( 1- \frac{\rho_+}{\rho_-} \right) g\eta .\label{weakly nonlinear press trans} \ee
Since $p$ is independent of the depth, this determines $\eta$ from the dynamic pressure on the top and bottom boundaries. 
In the case where the upper fluid has infinite extent, i.e., $H = +\infty$, we then have 
$p_{+} = 0$, and hence the dynamic pressure at the bed $p_{-}$ is proportional to $\eta$.  
This recovers the hydrodynamic approximation presented in Section \ref{sec:ww-hydrostatic}.

Moreover, the first equation in \eqref{E:euler4-scaled} in the moving coordinate becomes 
\be 
(U-c) u_x + v U_y = - p_x \qquad \textrm{in } \Omega,
\label{horizontal momentum 0} \ee
and \eqref{scaled kinematic cond} becomes 
\be v_\pm = (U_\pm - c) \eta_x \qquad \textrm{on } \{ y = h_0 \}. \label{two fluid lin traveling kinematic} \ee
We then use the last equation in \eqref{E:euler4-scaled} to rewrite \eqref{scaled kinematic cond} as 
\[ -(U-c) v_y + v U_y + p_x = 0 \qquad \textrm{in } \Omega,\]
which we rearrange as
\[  \rho \partial_y \left( \frac{v}{U-c} \right) = \rho \frac{p_x}{(U-c)^2} .\]
Here we have multiplied both sides by $\rho$ so that the jump conditions on the boundary are easy to interpret.  
Ultimately we integrate this and use \eqref{two fluid lin traveling kinematic} and \eqref{weakly nonlinear press trans}
to find the generalized Burns condition:
\be \int_{h_0}^{H} \frac{1}{(U_+-c)^2} \rho_+ \, dy + \int_0^{h_0} \frac{1}{(U_- -c)^2} \rho_- \, dy = 1. \label{generalized Burns} \ee
Here we tacitly assume that $\eta_x$ does not vanish identically 
and, in the case that $H = +\infty$, in addition, $1/(U_+ - c)$ is in $L^2$.  
When $\rho_+ = 0$, that is, the upper fluid region is a vacuum, 
\eqref{generalized Burns} reduces to \eqref{E:Burns}.

\subsection*{Acknowledgements}
RMC is supported in part by the Simons Foundation under Grant 354996 and the Central Research Development Fund No. 04.13205.30205 from University of Pittsburgh. VMH is supported by the National Science Foundation under CAREER DMS-1352597, an Alfred P. Sloan research fellowship, an Arnold O. Beckman research award RB14100 and a Beckman fellowship of the Center for Advanced Study at the University of Illinois at Urbana-Champaign. SW is supported in part by the National Science Foundation through grant DMS-1514910.

\bibliographystyle{amsalpha}
\bibliography{transferBib}

\def\cprime{$'$}
\providecommand{\bysame}{\leavevmode\hbox to3em{\hrulefill}\thinspace}
\providecommand{\MR}{\relax\ifhmode\unskip\space\fi MR }
\providecommand{\MRhref}[2]{%
  \href{http://www.ams.org/mathscinet-getitem?mr=#1}{#2}
}
\providecommand{\href}[2]{#2}
\begin{thebibliography}{OVDH12}

\bibitem[Ami84]{amick1984semilinear}
Charles~J. Amick, \emph{Semilinear elliptic eigenvalue problems on an infinite
  strip with an application to stratified fluids}, Ann. Scuola Norm. Sup. Pisa
  Cl. Sci. (4) \textbf{11} (1984), no.~3, 441--499. \MR{MR785621 (86i:35042)}

\bibitem[AT81]{AT1981solitary}
C.~J. Amick and J.~F. Toland, \emph{On solitary water-waves of finite
  amplitude}, Arch. Rational Mech. Anal. \textbf{76} (1981), no.~1, 9--95.
  \MR{629699 (83b:76017)}

\bibitem[AT86]{amick1986global}
C.~J. Amick and R.~E.~L. Turner, \emph{A global theory of internal solitary
  waves in two-fluid systems}, Trans. Amer. Math. Soc. \textbf{298} (1986),
  no.~2, 431--484. \MR{MR860375 (87m:35210)}

\bibitem[BD87]{bishop1987measuring}
Craig~T Bishop and Mark~A Donelan, \emph{Measuring waves with pressure
  transducers}, Coastal Engineering \textbf{11} (1987), no.~4, 309--328.

\bibitem[Bea77]{Beale1977solitary}
J.~Thomas Beale, \emph{The existence of solitary water waves}, Comm. Pure Appl.
  Math. \textbf{30} (1977), no.~4, 373--389. \MR{0445136 (56 \#3480)}

\bibitem[Bur53]{Burns}
J.~C. Burns, \emph{Long waves in running water}, Proc. Cambridge Philos. Soc.
  \textbf{49} (1953), 695--706, With an appendix by M. J. Lighthill.
  \MR{0057668 (15,261a)}

\bibitem[BWSZ15]{buhler2015wind}
Oliver B{\"u}hler, Samuel Walsh, Jalal Shatah, and Chongchun Zeng, \emph{On the
  wind generation of water waves}, Preprint (arXiv:1505.02032) (2015).

\bibitem[CC13]{clamond2013recovery}
Didier Clamond and Adrian Constantin, \emph{Recovery of steady periodic wave
  profiles from pressure measurements at the bed}, Journal of Fluid Mechanics
  \textbf{714} (2013), 463--475.

\bibitem[CEW07]{CEW2007symm}
Adrian Constantin, Mats Ehrnstr{\"o}m, and Erik Wahl{\'e}n, \emph{Symmetry of
  steady periodic gravity water waves with vorticity}, Duke Math. J.
  \textbf{140} (2007), no.~3, 591--603. \MR{2362244 (2009c:35359)}

\bibitem[CJ06]{constantin2006tsunami}
Adrian Constantin and Robin~Stanley Johnson, \emph{Modelling tsunamis}, J.
  Phys. A \textbf{39} (2006), no.~14, L215--L217. \MR{2219991 (2006k:86006)}

\bibitem[CJ08a]{constantin2008nondimen}
\bysame, \emph{On the non-dimensionalisation, scaling and resulting
  interpretation of the classical governing equations for water waves}, J.
  Nonlinear Math. Phys. \textbf{15} (2008), no.~suppl. 2, 58--73. \MR{2434725
  (2009h:76020)}

\bibitem[CJ08b]{constantin2008propagation}
\bysame, \emph{Propagation of very long water waves, with vorticity, over
  variable depth, with applications to tsunamis}, Fluid Dynam. Res. \textbf{40}
  (2008), no.~3, 175--211. \MR{2369543 (2009b:76009)}

\bibitem[CKS15]{constantin2015approximations}
Adrian Constantin, Konstantinos Kalimeris, and Otmar Scherzer,
  \emph{Approximations of steady periodic water waves in flows with constant
  vorticity}, Nonlinear Analysis: Real World Applications \textbf{25} (2015),
  276--306.

\bibitem[Cla13]{clamond2013new}
Didier Clamond, \emph{New exact relations for easy recovery of steady wave
  profiles from bottom pressure measurements}, Journal of Fluid Mechanics
  \textbf{726} (2013), 547--558.

\bibitem[Con12]{Constantin2012}
Adrian Constantin, \emph{On the recovery of solitary wave profiles from
  pressure measurements}, J. Fluid Mech. \textbf{699} (2012), 376--384.
  \MR{2923684}

\bibitem[CS04]{CS2004}
Adrian Constantin and Walter Strauss, \emph{Exact steady periodic water waves
  with vorticity}, Comm. Pure Appl. Math. \textbf{57} (2004), no.~4, 481--527.
  \MR{2027299 (2004i:76018)}

\bibitem[CS10]{constantin2010pressure}
\bysame, \emph{Pressure beneath a {S}tokes wave}, Comm. Pure Appl. Math.
  \textbf{63} (2010), no.~4, 533--557. \MR{2604871 (2011b:76017)}

\bibitem[CS11]{constantin2011discontinuous}
\bysame, \emph{Periodic traveling gravity water waves with discontinuous
  vorticity}, Arch. Ration. Mech. Anal. \textbf{202} (2011), no.~1, 133--175.
  \MR{2835865}

\bibitem[CSV14]{CSV2014}
Adrian Constantin, Walter Strauss, and Eugen Varvaruca, \emph{Global
  bifurcation of steady gravity water waves with critical layers}, 2014.

\bibitem[CW15a]{ChenWalsh2015pressure}
Robin~Ming Chen and Samuel Walsh, \emph{Reconstruction of stratified steady
  water waves from pressure readings}, Preprint, arXiv:1502.07775 (2015).

\bibitem[CW15b]{chen2015dependence}
\bysame, \emph{Continuous dependence on the density for stratified steady water
  wavesependence on the density for stratified steady water waves}, Arch.
  Rational Mech. Anal. (doi:10.1007/s00205-015-0906-6, 2015), 1--52.

\bibitem[DJ37]{dubreil1937theoremes}
ML~Dubreil-Jacotin, \emph{Sur les theoremes d'existence relatifs aux ondes
  permanentes periodiques a deux dimensions dans les liquides heterogenes}, J.
  Math. Pures Appl. \textbf{16} (1937), no.~9, 43--67.

\bibitem[DOV12]{deconinck2012relating}
B.~Deconinck, K.~L. Oliveras, and V.~Vasan, \emph{Relating the bottom pressure
  and the surface elevation in the water wave problem}, J. Nonlinear Math.
  Phys. \textbf{19} (2012), no.~suppl. 1, 1240014, 11. \MR{2999408}

\bibitem[DSP88]{da1988steep}
A.~F. Da~Silva and D.~H. Peregrine, \emph{Steep, steady surface waves on water
  of finite depth with constant vorticity}, Journal of Fluid Mechanics
  \textbf{195} (1988), 281--302.

\bibitem[EEW11]{EEW2011critical}
Mats Ehrnstr{\"o}m, Joachim Escher, and Erik Wahl{\'e}n, \emph{Steady water
  waves with multiple critical layers}, SIAM J. Math. Anal. \textbf{43} (2011),
  no.~3, 1436--1456. \MR{2821590 (2012f:76019)}

\bibitem[ES08]{ES2008}
Joachim Escher and Torsten Schlurmann, \emph{On the recovery of the free
  surface from the pressure within periodic traveling water waves}, J.
  Nonlinear Math. Phys. \textbf{15} (2008), no.~suppl. 2, 50--57. \MR{2434724
  (2009h:76022)}

\bibitem[FH54]{FH1954}
K.~O. Friedrichs and D.~H. Hyers, \emph{The existence of solitary waves}, Comm.
  Pure Appl. Math. \textbf{7} (1954), 517--550. \MR{0065317 (16,413f)}

\bibitem[Gro04]{Groves2004survey}
Mark~D. Groves, \emph{Steady water waves}, J. Nonlinear Math. Phys. \textbf{11}
  (2004), no.~4, 435--460. \MR{2097656 (2006a:76014)}

\bibitem[GW08]{GW2008solitary}
Mark~D. Groves and Erik Wahl{\'e}n, \emph{Small-amplitude {S}tokes and solitary
  gravity water waves with an arbitrary distribution of vorticity}, Phys. D
  \textbf{237} (2008), no.~10-12, 1530--1538. \MR{2454604 (2011a:37122)}

\bibitem[Hen13]{henry2013pressure}
David Henry, \emph{On the pressure transfer function for solitary water waves
  on the pressure transfer function for solitary water waves with vorticity},
  Math. Ann. \textbf{357} (2013), no.~1, 23--30.

\bibitem[HL08]{HL2008}
Vera~Mikyoung Hur and Zhiwu Lin, \emph{Unstable surface waves in running
  water}, Comm. Math. Phys. \textbf{282} (2008), no.~3, 733--796. \MR{2426143
  (2009d:76023)}

\bibitem[HL15]{HL2015pressure}
Vera~Mikyoung Hur and Michael~R. Livesay, \emph{On the recovery of traveling
  water waves with vorticity from the pressure at the bed}, Preprint (2015),
  arXiv:1510.02396.

\bibitem[Hur06]{Hur2006Stokes1}
Vera~Mikyoung Hur, \emph{Global bifurcation theory of deep-water waves with
  vorticity}, SIAM J. Math. Anal. \textbf{37} (2006), no.~5, 1482--1521
  (electronic). \MR{2215274 (2007e:76025)}

\bibitem[Hur07]{Hur2007symm}
\bysame, \emph{Symmetry of steady periodic water waves with vorticity}, Philos.
  Trans. R. Soc. Lond. Ser. A Math. Phys. Eng. Sci. \textbf{365} (2007),
  no.~1858, 2203--2214. \MR{2329142 (2008i:76025)}

\bibitem[Hur08]{Hur2008solitary}
\bysame, \emph{Exact solitary water waves with vorticity}, Arch. Ration. Mech.
  Anal. \textbf{188} (2008), no.~2, 213--244. \MR{2385741 (2010b:76019)}

\bibitem[Hur11]{Hur2011Stokes2}
\bysame, \emph{Stokes waves with vorticity}, J. Anal. Math. \textbf{113}
  (2011), 331--386. \MR{2788362 (2012e:76022)}

\bibitem[Hur12]{Hur2012regularity}
\bysame, \emph{Analyticity of rotational flows beneath solitary water waves},
  Int. Math. Res. Not. IMRN (2012), no.~11, 2550--2570. \MR{2926989}

\bibitem[Jam01]{james2001internal}
Guillaume James, \emph{Internal travelling waves in the limit of a
  discontinuously stratified fluid}, Arch. Ration. Mech. Anal. \textbf{160}
  (2001), no.~1, 41--90. \MR{1864121 (2002h:76031)}

\bibitem[Jan04]{janssen2004interaction}
P.~Janssen, \emph{The {I}nteraction of {O}cean {W}aves and {W}ind}, Cambridge
  University Press, 2004.

\bibitem[Joh97]{Johnson1997}
R.~S. Johnson, \emph{A modern introduction to the mathematical theory of water
  waves}, Cambridge Texts in Applied Mathematics, Cambridge University Press,
  Cambridge, 1997. \MR{1629555 (99m:76017)}

\bibitem[Kar12]{Pete2012dispersion}
Paschalis Karageorgis, \emph{Dispersion relation for water waves with
  non-constant vorticity}, Eur. J. Mech. B Fluids \textbf{34} (2012), 7--12.
  \MR{2927940}

\bibitem[KC94]{kuo1994transfer}
Yi-Yu Kuo and Yung-Fang Chiu, \emph{Transfer function between wave height and
  wave pressure for progressive waves}, Coastal Engineering \textbf{23} (1994),
  no.~1, 81--93.

\bibitem[Kin65]{kinsman1965wind}
Blair Kinsman, \emph{Wind {W}aves: their g{G}eneration and {P}ropagation on the
  {O}cean {S}urface}, Courier Corporation, 1965.

\bibitem[KS08a]{ko2008vorticity}
Joy Ko and Walter Strauss, \emph{Effect of vorticity on steady water waves}, J.
  Fluid Mech. \textbf{608} (2008), 197--215. \MR{2439751 (2009f:76023)}

\bibitem[KS08b]{ko2008large}
\bysame, \emph{Large-amplitude steady rotational water waves}, Eur. J. Mech. B
  Fluids \textbf{27} (2008), no.~2, 96--109. \MR{2389493 (2009b:76014)}

\bibitem[LF97]{lankers1997fast}
K.~Lankers and G.~Friesecke, \emph{Fast, large-amplitude solitary waves in the
  2d euler equations for stratified fluids}, Nonlinear Anal. \textbf{29}
  (1997), no.~9, 1061--1078.

\bibitem[OS01]{OS2001book}
Hisashi Okamoto and Mayumi Sh{\=o}ji, \emph{The {M}athematical {T}heory of
  {P}ermanent {P}rogressive {W}ater-{W}aves}, Advanced Series in Nonlinear
  Dynamics, vol.~20, World Scientific Publishing Co., Inc., River Edge, NJ,
  2001. \MR{1869386 (2003h:76016)}

\bibitem[OVDH12]{OVDH2012}
K.~L. Oliveras, V.~Vasan, B.~Deconinck, and D.~Henderson, \emph{Recovering the
  water-wave profile from pressure measurements}, SIAM J. Appl. Math.
  \textbf{72} (2012), no.~3, 897--918. \MR{2968755}

\bibitem[PT02]{PlotnikovToland2002}
Pavel~I. Plotnikov and John~F. Toland, \emph{The {F}ourier coefficients of
  {S}tokes' waves}, Nonlinear problems in mathematical physics and related
  topics, {I}, Int. Math. Ser. (N. Y.), vol.~1, Kluwer/Plenum, New York, 2002,
  pp.~303--315. \MR{1970618 (2004c:76023)}

\bibitem[Sch03]{schiereck2003introduction}
Gerrit~J Schiereck, \emph{Introduction to {B}ed, {B}ank and {S}hore
  {P}rotection}, CRC Press, 2003.

\bibitem[Sto47]{Stokes}
G.~G. Stokes, \emph{On the theory of oscillatory waves}, Trans. Camb. Phil.
  Soc. \textbf{8} (1847), 441--455.

\bibitem[Str10]{Strauss2010}
Walter~A. Strauss, \emph{Steady water waves}, Bull. Amer. Math. Soc. (N.S.)
  \textbf{47} (2010), no.~4, 671--694. \MR{2721042 (2012b:76022)}

\bibitem[SZ11]{shatah2011interface}
Jalal Shatah and Chongchun Zeng, \emph{Local well-posedness for fluid interface
  problems}, Arch. Ration. Mech. Anal. \textbf{199} (2011), no.~2, 653--705.
  \MR{2763036 (2012c:35336)}

\bibitem[TK63]{terkrikorov1963theorie}
A.~M. Ter-Krikorov, \emph{Th{\'e}orie exacte des ondes longues stationnaires
  dans un liquide h{\'e}t{\'e}rog{\`e}ne}, J. M{\'e}canique \textbf{2} (1963),
  351--376. \MR{MR0160400 (28 \#3613)}

\bibitem[Tol96]{Toland1996Stokes}
J.~F. Toland, \emph{Stokes waves}, Topol. Methods Nonlinear Anal. \textbf{7}
  (1996), no.~1, 1--48. \MR{1422004 (97j:35130)}

\bibitem[Tur81]{turner1981internal}
R.~E.~L. Turner, \emph{Internal waves in fluids with rapidly varying density},
  Ann. Scuola Norm. Sup. Pisa Cl. Sci. (4) \textbf{8} (1981), no.~4, 513--573.
  \MR{MR656000 (83j:76027)}

\bibitem[Tur84]{turner1984variational}
\bysame, \emph{A variational approach to surface solitary waves}, J.
  Differential Equations \textbf{55} (1984), no.~3, 401--438.

\bibitem[Var09]{Varvaruca2009extreme}
Eugen Varvaruca, \emph{On the existence of extreme waves and the {S}tokes
  conjecture with vorticity}, J. Differential Equations \textbf{246} (2009),
  no.~10, 4043--4076. \MR{2514735 (2010i:35450)}

\bibitem[VO14]{vasan2014pressure}
Vishal Vasan and Katie Oliveras, \emph{Pressure beneath a traveling wave with
  constant vorticity}, Discrete Contin. Dyn. Syst \textbf{34} (2014),
  3219--3239.

\bibitem[VW12]{VW2012extreme}
Eugen Varvaruca and Georg~S. Weiss, \emph{The {S}tokes conjecture for waves
  with vorticity}, Ann. Inst. H. Poincar\'e Anal. Non Lin\'eaire \textbf{29}
  (2012), no.~6, 861--885. \MR{2995099}

\bibitem[Wah09]{Wahlen2009critical}
Erik Wahl{\'e}n, \emph{Steady water waves with a critical layer}, J.
  Differential Equations \textbf{246} (2009), no.~6, 2468--2483. \MR{2498849
  (2010i:76026)}

\bibitem[Wal09]{walsh2009stratified}
Samuel Walsh, \emph{Stratified and steady periodic water waves}, SIAM J. Math.
  Anal. \textbf{41} (2009), no.~3, 1054--1105.

\bibitem[WBS13]{walsh2013wind}
Samuel Walsh, Oliver B{{\"u}}hler, and Jalal Shatah, \emph{Steady water waves
  in the presence of wind}, SIAM J. Math. Anal. \textbf{45} (2013), no.~4,
  2182--2227. \MR{3073648}

\bibitem[Whe13]{Wheeler2013solitary}
Miles~H. Wheeler, \emph{Large-amplitude solitary water waves with vorticity},
  SIAM J. Math. Anal. \textbf{45} (2013), no.~5, 2937--2994. \MR{3106477}

\bibitem[Whe15]{wheeler2015pressure}
\bysame, \emph{Solitary water waves of large amplitude generated by surface
  pressure}, Arch. Ration. Mech. Anal. \textbf{218} (2015), no.~2, 1131--1187.
  \MR{3375547}

\end{thebibliography}

\end{document}